\documentclass[a4paper,11pt]{amsart} 
\usepackage{amssymb,amsfonts,amsxtra, mathrsfs}
\usepackage[all]{xy}
\usepackage{fullpage}
\newtheorem{theorem}{Theorem}[section]
\newtheorem{cor}[theorem]{Corollary}

\newtheorem{prop}[theorem]{Proposition}

\theoremstyle{definition}

\newtheorem*{notation}{Notation}

\newtheorem{defi}[theorem]{Definition}
\newtheorem{rem}[theorem]{Remark} 

\numberwithin{equation}{section} 
\DeclareMathOperator{\Iso}{Iso}

\DeclareMathOperator{\Det}{Det}
\DeclareMathOperator{\Ori}{Or}
\DeclareMathOperator{\ori}{or}
\DeclareMathOperator{\st}{st}

\def\O{\mathscr O}
\def\F{\mathscr F}
\def\A{\mathscr A}
\def\M{\mathscr M}
\def\Mbar{\overline{\mathscr M}}
\def\T{\mathscr T}

\def\ra{\rightarrow} 
\def\R{\mathbb{R}}
\def\Z{\mathbb{Z}}
\def\s{\mathfrak{s}}

\newcommand{\noproof}{\begin{flushright} \ensuremath{\square} \end{flushright}}
\newcommand{\abs}[1]{{\lvert#1\rvert}}
\newcommand{\Ass}{\mathscr{A}\!\mathit{ss}}
\newcommand{\Comm}{\mathscr{C}\!\mathit{omm}}
\newcommand{\Lie}{\mathscr{L}\!\mathit{ie}}
\newcommand{\Out}{\operatorname{Out}}
\newcommand{\Rib}{\operatorname{Rib}}

\thanks{This work was completed during the first author's stay at IHES and he wishes to express his gratitude to this institution
for excellent working conditions.}
\begin{document}
\title[Koszul and Verdier duality...]{Graph homology: Koszul and
Verdier duality}
\author{A.~Lazarev \and A.~A.Voronov}
\address{IH\'ES, Le Bois-Marie, 35 route de Chartres, F-91440, Bures-sur-Yvette,
France.}
\email{lazarev@ihes.fr}
\address{School of Mathematics\\University of Minnesota\\206 Church
St. SE\\Minneapolis, MN 55455\\USA}
\email{voronov@umn.edu}

\begin{abstract}
  We show that Verdier duality for certain sheaves on the moduli
  spaces of graphs associated to Koszul operads corresponds to Koszul
  duality of operads. This in particular gives a conceptual
  explanation of the appearance of graph cohomology of both the
  commutative and Lie types in computations of the cohomology of the
  outer automorphism group of a free group. Another consequence is an
  explicit computation of dualizing sheaves on spaces of metric
  graphs, thus characterizing to which extent these spaces are
  different from oriented orbifolds. We also provide a relation
  between the cohomology of the space of metric ribbon graphs, known
  to be homotopy equivalent to the moduli space of Riemann surfaces,
  and the cohomology of a certain sheaf on the space of usual metric
  graphs.
\end{abstract}

\maketitle

\section*{Introduction}

The popularity of graph homology owes largely to the fact that the
cohomology of two important spaces in mathematics, the classifying
space $Y_n$ of the outer automorphism group of the free group on $n$
generators and the (decorated) moduli space $\M_{g,n}$ of Riemann
surfaces of genus $g$ with $n$ punctures, even though generally
intractable, may be computed via a deceptively simple combinatorial
construction, called \emph{graph homology}, see M.~Culler and K.~Vogtmann
\cite{CulV} and R.~C. Penner \cite{penner:bull}. These results,
combined with further study of graph homology by M.~Kontsevich
\cite{Kon}, rendered the following identifications:
\begin{eqnarray*}
  H_\bullet (Y_n, k) & \cong & {H}^{3n - 4 - \bullet}_{\Gamma \Lie}
  (n),\\
H^\bullet_c (Y_n, k) & \cong &
  \widetilde{H}^{\bullet}_{\Gamma \Comm} (n),\\
H_\bullet
  (\M_{g,n}, k) & \cong & {H}^{6g+3n-7 - \bullet}_{\Gamma \Ass} (g,n),\\
H^\bullet_c (\M_{g,n}, k) & \cong &
\widetilde{H}^{\bullet}_{\Gamma \Ass} (g,n),
\end{eqnarray*}
where $k$ is a coefficient field of characteristic zero, $H^\bullet_c$
denotes cohomology with compact supports, and in the right-hand side,
we have graph cohomology of various flavors, Lie, commutative, and
associative, with trivial or twisted coefficients.

The appearance of Koszul-dual operads in the right-hand side as
corresponding to the homology vs. cohomology with compact supports in the
left-hand side is quite suggestive: it hints on a relationship between
some kind of Poincar\'e duality for spaces and Koszul duality for
operads.

In this paper we show that this relationship indeed takes place and in
fact prove more general results, Theorems \ref{main2} and \ref{main4},
which show that up to an orientation twist, Verdier duality on the
moduli space of graphs transfers a certain constructible sheaf
corresponding to an operad $\O$ to the sheaf corresponding to the
dg-dual operad $D \O$, which is quasi-isomorphic to the Koszul-dual
operad $\O^!$, if $\O$ happens to be Koszul. The idea of a
relationship between the two dualities originates from the paper
\cite{GK} by Ginzburg and Kapranov, who noticed that Verdier duality
for sheaves on buildings (spaces of metric trees) provided a
sheaf-theoretic interpretation of Koszul duality for operads.
Koszulity was thus interpreted in terms of the vanishing of higher
cohomology for corresponding sheaves, while in our paper it translates
into a duality statement between highly nontrivial cohomology groups
of spaces of metric graphs.

The moduli spaces we consider are non-compact. It seems likely that
similar results hold for certain compactifications of our moduli
spaces and we intend to return to this in the future.

\begin{notation}
  Throughout this paper we work with vector spaces, graded vector
  spaces, and dg-vector spaces or complexes --- all finite-dimensional
  in each degree and bounded, over a ground field $k$, which is
  assumed to be of characteristic zero with the exception of
  Section~\ref{verdier}. We consider chain complexes $V_\bullet =
  \bigoplus_{i \in \Z} V_i$ with a differential $d: V_i \to V_{i-1}$
  and cochain complexes $V^\bullet = \bigoplus_{i \in \Z} V^i$ with a
  differential $d: V^i \to V^{i+1}$.

  The \emph{$($degree$)$ shift} $V[1]$ of a complex $V$ has components
  $(V[1])_i = V_{i+1}$ in the category of chain complexes and
  $(V[1])^i = V^{i+1}$ in the category of cochain complexes. For chain
  complexes the degree shift is also known as desuspension.

  The functor $V \mapsto V^*$ of taking the linear dual acts within
  each of the two categories:
\begin{eqnarray*}
V^*_i & = (V_{-i})^*, & d^*: V^*_i \to V^*_{i-1},\\
(V^*)^i & = (V^{-i})^*, & d^*: (V^*)^i \to (V^*)^{i+1},
\end{eqnarray*}
   while another functor, $V \mapsto V^\vee$, takes the category of
   chain complexes to that of cochain ones:
\[
(V^\vee)^i  = (V_{i})^*,  d^\vee: (V^\vee)^i \to (V^\vee)^{i+1}.
\]
Note that $(V[1])^* \cong V^*[-1]$ and $(V[1])^\vee \cong V^\vee[1]$.
The double dual $V^{**}$ of a chain complex $V$ is naturally
isomorphic to $V$, while $V^{* \vee} \cong V^{\vee *}$ and the functor
$V\mapsto V^{\vee *}$ is an equivalence of categories of chain and
cochain complexes. Clearly $(V^{\vee *})^i \cong V_{-i}$.

An ungraded vector space $V$ could be assumed to lie in degree $0$,
and it will be clear from the context whether this (trivial) grading
is considered homological or cohomological. If $\dim V=n$ we will call
the \emph{determinant} of $V$ the one-dimensional graded vector space
$\Det(V) = S^n(V[-1]) = \Lambda^n (V)[-n]$, concentrated in degree
$n$. Note that $\Det(V)^*[-2n] \cong  \Det (V^*)$.
We will use negative powers of one-dimensional graded vector spaces
for the corresponding positive tensor powers of their $*$-duals, so
that
\[
\Det^{-p} (V) = \left( (\Det V)^* \right)^{\otimes p}.
\]

For a finite collection $\{ V_\alpha \; | \; \alpha \in I\}$ of
finite-dimensional vector spaces, we have a natural identification
\[
\bigotimes_{\alpha \in I} V_\alpha [-1] \cong \Det (I) \otimes
\bigotimes_{\alpha \in I} V_\alpha.
\]

If $S$ is a finite set, let $\Det (S) := \Det (k^S)$. Since there is a canonical isomorphism $(k^S)^* \cong
k^S$, we have $\Det (S)^*[-2\abs{S}] \cong \Det
(S)$. Note also that $\Det^2 (S) \cong k [-2\abs{S}]$.

For a simplex $\sigma$, the symbol $\Det(\sigma)$ will denote the determinant of
the set of vertices of $\sigma$.  When the ground field $k = \R$, a
choice of a nonzero element in $\Det(\sigma)$ up to a positive real
factor is equivalent to providing $\sigma$ with an orientation in the
usual sense.
\end{notation}

\section{Verdier duality for simplicial complexes}
\label{verdier}

In this section we formulate and prove certain results on Verdier
duality for sheaves on simplicial complexes. These results, in a
slightly different situation of spaces stratified into cells, were
stated in \cite{GK}.

\begin{defi} Let $X$ be a finite simplicial complex.  A sheaf of
  dg-vector spaces over a ground field $k$ on $X$ is called
  \emph{constructible}, if its restriction to each face of $X$ is a
  constant sheaf whose stalk is a dg-vector space.
\end{defi}

\begin{rem}
  Ginzburg and Kapranov use the term ``combinatorial sheaf.'' We
  follow the more conventional terminology adopted in, e.g. \cite{KS}.
\end{rem}
Any simplicial complex $X$ admits an open covering $U_\sigma$ where
$\sigma$ runs through the faces of $X$; namely $U_\sigma$ is the open
star of $\sigma$, the union of the interiors of those faces of $X$ which
contain $\sigma$. Any sheaf determines a contravariant functor from
the poset $\{U_\sigma\}$ into the category of dg-vector spaces.
Conversely, let $\F$ be a constructible sheaf on $X$. Let $x\in X$ and
consider the face $\sigma$ of smallest dimension containing $x$. Then
the set of sections of $\F$ over any sufficiently small neighborhood
of $X$ will coincide with $\Gamma(X,U_\sigma)$. Therefore, $\F$ is
completely determined by the corresponding functor.

Consider the category whose objects are the simplices of $X$ and the
morphisms are inclusions of faces. We will call a \emph{coefficient
  system} on $X$ any covariant functor from this category to the
category of dg-vector spaces.
\begin{prop}
  There is a one-to-one correspondence between constructible sheaves
  and coefficient systems on a simplicial complex $X$.
\end{prop}
\begin{proof}
Indeed, it suffices to note that the category of faces of $X$ is opposite to the category of open stars of $X$.
\end{proof}

\begin{rem}
\label{acyclicity}
The cohomology of a constructible sheaf could be computed using the
\v{C}ech complex of the covering $\{U_\sigma\}$, as follows from
Kashiwara-Schapira \cite{KS}, Proposition 8.1.4. Cohomology in this
paper will always mean hypercohomology.
\end{rem}

\begin{prop}
\label{GK1}
Let $\F$ be a constructible sheaf and $\{\F_\sigma\}$ be the
corresponding coefficient system on $X$.  Then the cohomology of a
constructible sheaf on $X$ coincides with the cohomology of the cochain
complex
\begin{equation}
\label{cech}
\bigoplus_{\tau}\F_\tau \otimes \Det(\tau)[1]
\end{equation}
on which the differential acts as the sum of the internal differential
on $\F$ and a map
\[
\F_\sigma
\otimes \Det(\sigma)[1] \mapsto\sum_{\begin{subarray}{c}\tau\supset\sigma\\
    \dim\tau=\dim\sigma+1\end{subarray}} \F_{\tau} \otimes
\Det(\tau)[1],
\]
where the last map is induced by inclusions $\sigma\hookrightarrow
\tau$.
\end{prop}
\begin{proof}
  According to Remark~\ref{acyclicity}, the cohomology of $\F$ could
  be computed using the \v{C}ech (bi)complex of the covering of
  $\{U_\sigma\}$ of $X$. A simple inspection shows that this complex
  is isomorphic to the complex (\ref{cech}).
\end{proof}
We will now discuss  Verdier duality in the simplicial
context. Recall that for a sheaf $\F$, considered as an object of the
derived category of sheaves on $X$, its \emph{Verdier dual} $D\F$ is
characterized by the property $R\Gamma(U,D\F)=[R\Gamma_c(U,\F)]^*$.
It is easy to see that for a constructible sheaf $\F$, its Verdier
dual complex $D\F$ will have constructible cohomology and therefore,
by \cite{KS}, Theorem 8.1.10, can be represented by a complex of
constructible sheaves.
\begin{prop}\label{GK2}
  Let $\F$ be a constructible sheaf on $X$ and $D\F$ is its Verdier
  dual. Then $D\F$ is represented by the constructible complex
  $\sigma\mapsto D\F_\sigma$, where $D\F_\sigma$ is the following
  cochain complex:
\[
\bigoplus_{\tau\supset\sigma} (\F_\tau\otimes \Det(\tau)[1])^*
\]
whose differential is the dual to that in \eqref{cech}.
\end{prop}
Note that under our grading convention for dual spaces, $\deg
(\Det (\tau)[1])^* = -\dim\tau$.
\begin{proof}
  Consider the open star $\st(\sigma)$ of the simplex $\sigma$.  We
  will denote by $i:\st(\sigma)\rightarrow \overline{\st(\sigma)}$ the
  inclusion of $\st(\sigma)$ into its closure. Then
  $i_!(\F|_{\st(\sigma)})$ is a constructible sheaf on the simplicial
  complex $\overline{\st(\sigma)}$. It follows that
  \[R\Gamma_c(\F,\st(\sigma))=R\Gamma_ci_!(\F|_{\st(\sigma)})=R\Gamma
  i_!(\F|_{\st(\sigma)}).\] Note that $\overline{st(\sigma)}$ is the
  union of all simplices containing $\sigma$. The sheaf
  $i_!(\F|_{\st(\sigma)})$ corresponds to the coefficients system on
  $\overline{st(\sigma)}$ so that
  $(i_!\F|_{\st(\sigma)})_\tau=\begin{cases}\F_\tau~ \text{if}~
    \tau\supset\sigma,\\ 0 ~\text{otherwise.}\end{cases}$ Now
  Corollary \ref{GK1} implies that $R\Gamma i_!(\F|_{\st(\sigma)})$ is
  represented by the complex
\[
\bigoplus_{\tau\supset\sigma} \F_\tau\otimes  \Det(\tau)[1],
\]
and the desired statement follows.
\end{proof}
\begin{rem}
  This result was formulated in the stratified setting in \cite{GK},
  Proposition 3.5.12 (b).
\end{rem}

\section{Equivariant Verdier duality}
\label{e-verdier}

In this section we generalize our theory to the case of
\emph{orbi-simplicial} complexes. We will not discuss orbi-simplicial
complexes in full generality, restricting ourselves to the case when
there exists a global group action. For the rest of the paper, the
ground field $k$ will have characteristic $0$.  Let $X$ be a
topological space and $G$ be a group acting properly discontinuously
on $X$.  That means that the stabilizer $G_x$ of every point $x\in X$
is finite and every point $x\in X$ has a neighborhood $U_x$ such that
$gU_x\bigcap U_x=\emptyset$ if $g\notin G_x$. Let $Y$ denote the space
of orbits $X/G$ and by $f:X\ra Y$ the projection map.  We now recall
some standard definitions and facts about equivariant sheaves, cf.
\cite{Tohoku} or a more modern reference \cite{Lunts}.
\begin{defi}
  A $G$-equivariant sheaf $\F$ on $X$ is a sheaf of $k$-vector spaces
  with a $G$-action. More precisely, for any $g\in G$ and any open set
  $U\subset X$ there is an isomorphism $g_U:\Gamma(U,\F)\rightarrow
  \Gamma(gU, \F)$ which is compatible with the restriction maps in the
  sense that for any open subsets $V\subset U$ in $X$ the following
  diagram is commutative:
\[
\xymatrix{\Gamma(U, \F)\ar^{g_U}[r]\ar[d]&\Gamma(gU, \F)\ar[d]\\\Gamma(V, \F)\ar^{g_V}[r]&\Gamma(gV, \F)}
\]
where the downward arrows are the restriction maps.  In addition we
require the following cocycle conditions:
\begin{itemize}
\item $1_U$ is the identity isomorphism for any open set $U$.\item
$h_{gU}\circ g_{U}=(h\circ g)_{U}$ for any $h,g\in G$ and any open set
$U\in X$.
\end{itemize}
\end{defi}
Note that $\Gamma(X, \F)$ has a $G$-action. We will denote by
$\Gamma^G(X, \F)$ the space of $G$-invariants: $\Gamma^G(X,
\F)=[\Gamma(X, \F)]^G$. A morphism $\F_1\ra \F_2$ between two
equivariant sheaves is an element in
$\Gamma^G(X,\mathscr{H}\!\mathit{om}(\F_1,\F_2))$. $G$-equivariant
sheaves on $X$ form an abelian category.  For any sheaf $ \F$ on $Y$
the sheaf $f^{-1} \F$ is naturally a $G$-equivariant sheaf on $X$.  The
direct image sheaf $f_* \F$ is a $G$-equivariant sheaf on $Y$ where $G$
is assumed to act trivially on $Y$.
\begin{defi} 
The $G$-equivariant direct image $f_*^G \F$ is the sheaf of
$G$-invariants of $f_* \F$ so that for $V\in Y$ we have $\Gamma(V,f_*^G
\F)=\Gamma(V,f_*\F)^G$.
\end{defi}
The functor $f^{-1}$ embeds the category of sheaves on $Y$ as a full
subcategory into the category of $G$-equivariant sheaves on
$X$. Moreover, $f_*^G\circ f^{-1}$ is isomorphic to the identity
functor on the category of sheaves on $Y$. Since the functor $f_*^G$
is exact these statements continue to hold on the level of derived
categories, cf. \cite{Lunts}, Theorem 8.6.1.

Now assume that $X$ is a finite-dimensional simplicial complex and
that $G$ acts simplicially, i.e. for any simplex $\sigma\in X$ and
$g\in G$ the image $g(\sigma)$ is another simplex of $X$ and
$g:\sigma\rightarrow g(\sigma)$ is an affine map. Our standing
assumptions on the action imply that the stabilizer of each simplex is
finite.  As a topological space $Y$ is glued from
\emph{orbi-simplices}, i.e. quotients of simplices by actions of
finite groups. One has one $n$-dimensional orbi-simplex of $Y$ for
each orbit of the action of $G$ on the set of $n$-simplices of $X$.
\begin{defi}
A sheaf $\F$ on $Y$ is called constructible if $f^{-1}\F$ is constructible on $X$.
\end{defi}
In other words a constructible sheaf is constant when restricted onto
each orbi-simplex. Just as in the non-equivariant situation, a
constructible sheaf $\F$ on $Y$ is equivalent to a coefficient system
on $Y$, i.e. a functor $\sigma\mapsto \F_\sigma$ from the poset of
orbi-simplices of $Y$ into $k$-vector spaces.  Then we have the
following (almost verbatim) analogue of Corollary \ref{GK1}.
\begin{prop}
\label{GKa1}
Let $\F$ be a constructible sheaf on $Y$ and $\{\F_\sigma\}$ be the
corresponding coefficient system on $Y$.  Then the cohomology of a
constructible sheaf on $X$ coincides with the cohomology of the
complex
\begin{equation}
\label{checheq}
\bigoplus_{\tau}\F_\tau \otimes \Det(\tau)[1].
\end{equation}
Here the direct sum is over the orbi-simplices of $Y$, and the
differential acts as in the non-equivariant situation.
\end{prop}
\begin{proof}
  According to the correspondence between equivariant sheaves on $X$
  and non-equivariant sheaves on $Y$ we have an isomorphism
  $R\Gamma^G(X, f^{-1}\F)\cong R\Gamma(Y,\F)$. Since the complex
  (\ref{checheq}) is just the complex of $G$-invariants of the
  \v{C}ech complex of $f^{-1}\F$ and the latter does compute
  $R\Gamma(X,\F)$ by Corollary \ref{GK1} the statement of our
  proposition follows.
\end{proof}
Similar arguments can be used to prove the following analogue of Proposition \ref{GK2}.
\begin{prop}
\label{GKa2}
Let $\F$ be a constructible sheaf on $Y$ and $D\F$ be its Verdier
dual. Then $D\F$ is represented by the constructible complex
$\sigma\mapsto D\F_\sigma$ where $D\F_\sigma$ is the following
complex:
\begin{equation}
\label{ver}
\bigoplus_{\tau\supset\sigma} (\F_\tau\otimes \Det(\tau)[1])^*.
\end{equation}
Here $\tau$ runs over the orbi-simplices of $Y$ having $\sigma$ as a
face, the grading convention and the formula for the differential are
the same as in the non-equivariant situation.
\end{prop}\noproof

\begin{rem}
\label{orbisimplicial}
We need a mild generalization of the above results. Suppose that
$X^\prime$ is a $G$-subcomplex of $X$ and $Y^\prime:=X^\prime/G\subset
Y=X/G$. We require that $G$ act properly discontinuously on
$X\setminus X^\prime$, but \emph{not necessarily} on the whole $X$. If
$\sigma$ is a simplex in $X\setminus X^\prime$ we will still refer to
its image in $Y$ as an \emph{orbi-simplex}.

Furthermore, let $\F$ be a sheaf on $Y$ for which $i^{-1}\F$ is a
constructible sheaf on $X$ and $\F|_{Y^\prime}=0$.

Then Propositions \ref{GKa1} and \ref{GKa2} continue to hold. In other
words formula (\ref{checheq}) still computes the cohomology of $\F$
and (\ref{ver}) represents the sheaf $D\F$.  To see that one only has
to note that
\begin{itemize}
\item the functor $i^{-1}$ restricted to the category of sheaves on
  $Y\setminus Y^\prime$ still embeds the derived category of sheaves
  on $Y\setminus Y^\prime$ into the derived category of sheaves on
  $X\setminus X^\prime$ and
\item for the inclusion $j:X\setminus X^\prime\hookrightarrow X$ the functor
  $j_!$, the direct image with compact support is exact.
\end{itemize}
\end{rem}

\section{Graph complexes and spaces of metric graphs}
\label{metric}

\subsection{Graph complexes}
\label{graphs}
A graph is specified by a set of vertices, a set of half-edges, and
(rather obvious) combinatorial relations between them, cf.\ for
example, \cite{GetK:modular} for precise definitions. One may also think of a
graph as an isomorphism class of a $1$-dimensional CW complex. We will
only consider connected, finite graphs whose vertices have valence
three or higher, i.e., for each vertex the number of incident
half-edges must be at least three. The sets of vertices and edges of a
graph $\Gamma$ will be denoted by $V(\Gamma)$ and $E(\Gamma)$,
respectively. The set of half-edges incident to a vertex $v\in
V(\Gamma)$ will be denoted by $H(v)$.

Let $\O$ be a cyclic operad in the category of chain complexes of
$k$-vector spaces.  For simplicity, we will assume that $O(1) = k$ and
$\O(n)$ is a finite-dimensional dg-vector space for each $n \ge 2$, as
this is the case for the standard examples of $\O = \Comm$, $\Ass$,
and $\Lie$. The more general case of an admissible operad, see
\cite{GK}[3.1.5], can also be treated by taking tensor products over
the associative algebra $K = \O(1)$, rather than the ground field $k$.

If $S$ is a set of $n+1$ elements, $n > 0$, one can define $\O((S))$
by using the coinvariants trick:
\[
\O((S)) := (\O(n) \times \Iso(S,[n]))_{S_{n+1}},
\]
where $\Iso (S,[n])$ is the set of bijections between $S$ and $[n] :=
\{0, 1, \dots, n\}$ and the symmetric group $S_{n+1}$ acts
diagonally. Recall the notion of an $\O$-graph complex \cite{CV}.

\begin{defi}
  An $\O$-\emph{decorated graph} or simply an $\O$-\emph{graph} is a
  graph $\Gamma$ together with a decoration which associates to any
  vertex $v$ of $\Gamma$ an element in $\O((H(v)))$.

  The \emph{space of $\O$-decorations on} $\Gamma$ is the
  chain complex
\[
\Gamma^{\O} = \bigotimes_{v\in V(\Gamma)}{\O}((H(v))).
\]

The \emph{orientation space} of a graph $\Gamma$ is the
one-dimensional graded vector space $\Ori(\Gamma):=\Det (E(\Gamma))
\otimes \Det^{-1} H_1(\Gamma)[\chi]$, concentrated in degree
$e(\Gamma) - 1$, where $e(\Gamma) = \abs{E(\Gamma)}$ and $\chi =
\chi(\Gamma)$ is the Euler characteristic of the graph $\Gamma$ as a
CW complex.  A \emph{twisted orientation space} of a graph $\Gamma$ is
the vector space $\Det(E(\Gamma))[1]$.  A \emph{(twisted) orientation}
on a graph $\Gamma$ is a choice of a nonzero element $\ori$ in
$\Ori(\Gamma)$ ($\Det(E(\Gamma)[1])$, respectively).
\end{defi}

\begin{rem}
  When $k = \R$, an orientation on a graph (up to a positive real
  factor) is equivalent to an ordering of its vertices and directing
  its edges (up to even permutation), cf.\
  \cite{GetK:modular,Thur,CV}.
\end{rem}

The following cyclic operads are of particular importance:
\begin{enumerate}
\item The commutative operad $\Comm(n)=k$ for $n>0$;
\item The associative operad $\Ass(n)=k[S_n]$ for $n>0$;
\item The Lie operad, whose $n$th space $\Lie(n)$ is the $k$-vector
  space spanned by all Lie monomials in $n$ variables containing each
  variable exactly once.
\end{enumerate}
The corresponding $\O$-graphs are called commutative, ribbon
and Lie graphs, respectively.
\begin{defi}
  The $\O$-\emph{graph complex} is the following complex of
  $k$-vector spaces:
\[
C_\bullet^{\Gamma \O} = \bigoplus_\Gamma \Gamma^\O \otimes
\Ori(\Gamma),
\]
where the summation runs over the isomorphism classes $\Gamma$ of
graphs. The grading comes from the internal grading on $\O$ and the
grading on $\Ori(\Gamma)$, which sits in degree $e(\Gamma) - 1$, so
that, provided that $\O$ is non-negatively graded, the graph complex
in general would end in degree $-\chi(\Gamma)$, corresponding to
graphs with one vertex.  The differential is the sum of the internal
differential coming from the operad $\O$ and the graph differential
$d:C_n^{\Gamma \O}\rightarrow C_{n-1}^{\Gamma \O}$ which acts as
follows:
\begin{equation}
\label{graph-d}
d(\Gamma \otimes \ori )=\sum_e \Gamma_e \otimes \ori_e,
\end{equation}
where $\Gamma$ is an $\O$-graph and $\ori \in \Ori(\Gamma)
\setminus \{0\}$ is an orientation on $\Gamma$. Here $\Gamma_e$ is the
graph obtained from $\Gamma$ by contracting an edge $e$ and the
summation is taken over all edges of $\Gamma$ which are not loops.

The orientation $\ori_e$ of $\Gamma_e$ is induced from the orientation
$\ori$ of $\Gamma$ in such a way that $\ori = e \wedge \ori_e$. The
$\O$-decoration on $\Gamma$ is defined as follows. Suppose that the
two endpoints $v_1$ and $v_2$ of the edge $e$ have valences $n_1$ and
$n_2$, respectively. Then the vertex obtained by coalescing $v_1$ and
$v_2$ is decorated by the element $v_1(\O)\circ v_2(\O)\in
\O(n_1+n_2-3)$, where $\O(n_1-1)\circ \O(n_2-1)\rightarrow
\O(n_1+n_2-3)$ is a structure map of the operad
$\O$.

Similarly, the \emph{twisted} graph complex
$\widetilde{C}_\bullet^{\Gamma \O}$ is formed by the isomorphism
classes of $\O$-graphs with twisted orientation; the grading and the
differential are defined like in the untwisted case, so that, for
example, terms corresponding to graphs $\Gamma$ with a single vertex
decorated by an element of $\O$ of degree zero would sit in degree
$-\chi(\Gamma)$.

The homology of the complexes ${C}_\bullet^{\Gamma \O}$ and
$\widetilde{C}_\bullet^{\Gamma \O}$ are denoted by
${H}_\bullet^{\Gamma \O}$ and $\widetilde{H}_\bullet^{\Gamma \O}$
respectively and called \emph{$\O$-graph homology}. The cohomology of
the $k$-dual cochain complexes ${C}^\bullet_{\Gamma \O} =
[{C}_\bullet^{\Gamma \O}]^\vee = \bigoplus_\Gamma (\Gamma^\O)^\vee
\otimes \Det (E(\Gamma)) \otimes \Det^{-1} (H_1(\Gamma)^*)[\chi]$ and
$\widetilde{C}^\bullet_{\Gamma \O} = [\widetilde{C}_\bullet^{\Gamma
  \O}]^\vee = \bigoplus_\Gamma (\Gamma^\O)^\vee \otimes \Det
(E(\Gamma))[1]$ are called \emph{$\O$-graph cohomology} and
\emph{twisted $\O$-graph cohomology,} respectively.
\end{defi}
\begin{rem}
  The graph complexes associated to two quasi-isomorphic operads are
  themselves quasi-isomorphic as a standard spectral sequence argument
  makes clear.
\end{rem}
For an integer $n>1$, we will consider graphs $\Gamma$ with
$H_1(\Gamma, \Z)$ being a free abelian group of rank $n$. These graphs
form a subcomplex ${C}_\bullet^{\Gamma \O}(n)$; clearly
${C}_\bullet^{\Gamma \O}\cong\bigoplus_n {C}_\bullet^{\Gamma \O}(n)$.
In what follows the number $n$ will be understood but not explicitly
mentioned. Note that ${C}_\bullet^{\Gamma \O}$ are complexes of
finite-dimensional vector spaces and have finite lengths.

Furthermore in the case $\O = \Ass$ an $\O$-graph -- a ribbon graph
-- has an additional invariant, the \emph{number $b > 0$ of boundary
components}, cf.\ for example, the survey \cite{HL}.  It is convenient
to introduce the \emph{genus} $g \ge 0$ of a ribbon graph by the
formula $g=1/2(n+1-b)$; then graphs with fixed $g$ and $b$ form a
subcomplex ${C}_\bullet^{\Gamma \Ass}(g,b)$ inside
${C}_\bullet^{\Gamma \Ass}$ and ${C}_\bullet^{\Gamma
\Ass} (n) \cong
\bigoplus_{\substack{g \ge 0, b \ge 1\\2-2g-b = 1 -n < 0}}
{C}_\bullet^{\Gamma \Ass}(g,b)$.

\subsection{Metric graphs}
We now introduce the moduli space of metric graphs, cf.\ \cite{CulV}.
A \emph{metric graph} is a graph $\Gamma$ together with a map
$l:E(\Gamma)\rightarrow \mathbb{R}_+$; the positive number $l(e)$ is
called the \emph{length} of the edge $e$.

A \emph{marking} on a metric graph is a homotopy equivalence $\vee_1^n
S^1\rightarrow \Gamma$ from an $n$-fold wedge of circles to $\Gamma$.
The set of markings on $\Gamma$ is acted on by the group of isometries
of $\Gamma$ and marked graphs belonging to the same orbit are called
equivalent.

We associate to any metric graph $\Gamma$ with $\sum_{e\in E(\Gamma)}
l(e) = 1$ a point in the simplex $\Delta^{e(\Gamma)-1}$ whose
barycentric coordinates are given by the lengths of the edges of
$\Gamma$. Collapsing edges corresponds to passing to the faces of the
corresponding simplex and after a suitable identification we obtain a
topological space $X_n$, the so-called \emph{Outer Space}. A point in
$X_n$ is a marked graph and the group of outer automorphisms of the
free group on $n$ generators, denoted by $\Out(F_n)$, acts on $X_n$ by
changing markings. The space $X_n$ is not a simplicial complex owing
to the existence of loops but adding formally the missing faces (ideal
simplices), one obtains a simplicial complex $\overline{X}_n$, on
which the group $\Out(F_n)$ continues to act, together with an
inclusion $i:Y_n\rightarrow\overline{Y}_n$.

The space $X_n$ is contractible and $\Out(F_n)$ acts on it properly
discontinuously, so $Y_n:=X_n/\Out(F_n)$ is rationally a classifying
space of $\Out(F_n)$. Each isomorphism class of graphs contributes an
orbi-simplex in $\overline{Y}_n$. The remaining simplices will be
referred to as \emph{ideal orbi-simplices}.

Note that $\overline{Y}_n :=Y_n/\Out(F_n)$ is \emph{not} a simplicial
orbi-complex since $\Out(F_n)$ does not act on the ideal simplices
properly discontinuosly. However this will not affect our results
since the sheaves we are interested in will vanish on those bad
simplices, cf.\ Remark \ref{orbisimplicial}.  We will now introduce
certain constructible sheaves on $\overline{\mathscr H}$ on
$\overline{Y}_n$.
\begin{defi}\
\begin{enumerate}
\item For an orbi-simplex $\sigma$ corresponding to a graph $\Gamma$,
  we set $\overline{{\mathscr H}}_\sigma=\Det^{-1} (H_1(\Gamma))
  [-n]$. If $\sigma$ is an ideal orbi-simplex we set
  $\overline{{\mathscr H}}_\sigma=0$. The complex
  $i^{-1}\overline{\mathscr{H}}$ on $Y_n$ will be denoted by
  ${\mathscr H}$.
\item Associated to a cyclic operad $\O$ is a constructible complex
  $\overline{\F}^{\O}$ on $\overline{Y}_n$ defined as follows. For an
  oriented simplex $\sigma\in \overline{Y}_n$ corresponding to a graph
  $\Gamma$, we set $\F^\O_\sigma$ to be the cochain complex $k$-dual
  to the chain complex $\Gamma^\O$ of $\O$-decorations on $\Gamma$:
\[
\F^\O_\sigma := (\Gamma^\O)^\vee.
\]
If $\sigma$ is an ideal simplex we set $\F^{\O}_\sigma=0$. If
$\sigma \subset \tau$ is a face of $\tau$, the corresponding morphism
$\F^\O_\sigma \to \F^\O_\tau$ is defined as the dual to the one
obtained by the operad composition of decorations along the edges in
the graph $\Gamma_\tau$ being contracted to obtain $\Gamma_\sigma$,
where $\Gamma_\sigma$ and $\Gamma_\tau$ are the graphs corresponding
to the simplices $\sigma$ and $\tau$, respectively. The complex of
sheaves $i^{-1}\overline{\F}^{\O}$ on $Y_n$ will be
denoted by ${\F}^{\O}$.
\end{enumerate}
\end{defi}
\begin{rem}\
\begin{itemize}
\item Clearly $i_!{\F}^{\O}\cong
  \overline{\F}^{\O}$; similarly $i_!{\mathscr
    H}\cong \overline{\mathscr H}$. Furthermore ${\mathscr H}$ is a
  locally free sheaf on $Y_n$.
\item For two quasi-isomorphic operads $\O$ and $\O'$, the complexes
  $\F_{\O}$ and $\F_{\O'}$ are quasi-isomorphic.
\end{itemize}
\end{rem}
\begin{theorem}
\label{main1}
There are canonical isomorphisms of graded $k$-vector spaces:
\begin{enumerate}
\item $\widetilde{H}^\bullet_{\Gamma\O}(n)\cong H^{\bullet}_c(Y_n,
  \F^{\O})$ for the twisted $\O$-graph cohomology;
\item ${H}^\bullet_{\Gamma\O}(n)\cong H^{\bullet}_c(Y_n,
  \F^{\O}\otimes \mathscr H)$ for the $\O$-graph cohomology.
\end{enumerate}
\end{theorem}

\begin{proof} We shall only prove part (1), the argument for (2) being
  virtually identical, as the standard orientation on a graph $\Gamma$
  differs from the twisted one by $\Det^{-1}(H_1(\Gamma))[-n]$. It
  suffices to prove the isomorphism
  $\widetilde{H}^\bullet_{\Gamma\O}(n) \cong
  H^{\bullet}(\overline{Y}_n, \overline{\F}^{\O})$.  Using Proposition
  \ref{GKa1} (or, rather, its modification as in
  Remark~\ref{orbisimplicial}), we see that the complex computing
  $H^\bullet(\overline{Y}_n,\overline{\F}^{\O})$ coincides with the
  complex $\widetilde{C}^\bullet_{\Gamma\O}(n)$, so the statement of
  the theorem follows.
\end{proof}

Given a cyclic operad $\O$ of chain complexes, let us
recall the construction of its \emph{dg-dual operad} $D\mathscr O$
from \cite{GK,GetK:cyclic,GetK:modular}. The $n$th component $D \O
(n)$ of the dg-dual cyclic operad $D \O$ for each $n > 0$ is defined
as the linear dual of the space of oriented, unrooted trees decorated
by a certain degree shift $\s \O[-1]$ with leaves labeled by $0, 1,
\dots, n$. More precisely,
\begin{equation}
\label{dg-dual}
D \O(n) := \bigoplus_{\text{unrooted $n+1$-trees T}}
(T^{\s\O[-1]})^*,
\end{equation}
where the summation runs over the isomorphism classes $T$ of
(unrooted) trees with vertices of valence at least three and $n+1$
leaves thought of as half-edges with free ends and labeled by numbers
$0, 1, \dots, n$. Here $\s\O$ is the \emph{cyclic-operad
suspension},
\cite{GetK:cyclic}:
\[
\s\O(n) := \Det^{-1}(k^{n+1}) [-2] \otimes \O(n),
\]
which results in
\[
\s\O((S)) = \Det^{-1}(S) [-2] \otimes \O((S))
\]
for a finite set $S$. Also, $T^{\s \O[-1]}$ is the space of $\s
\O[-1]$-decorations on $T$. Thus, if $\O$ happens to be concentrated
in degree zero, $D \O(n)$ will be a chain complex spanning degrees
$n-2$ through $0$. The differential on $D\O(n)$ is the sum of the
internal differential coming from the complex of $\O$-decorations on a
tree and the differential linear dual to the differential
\eqref{graph-d} restricted from graphs to trees.

\begin{rem}
  If $S$ is a set of $n+1$ elements, one can make precise sense out
  of $D \O(S)$ by considering trees whose leaves are labeled by the
  elements of $S$.
\end{rem}

\begin{theorem}
\label{main2}
There is a canonical isomorphism in the derived category of sheaves on
${Y}_n$:
\[
D\F^\O \cong \F^{D\O} \otimes {\mathscr H} [4-3n],
\]
where $D \F$ is the Verdier dual sheaf and $D\O$ the dg-dual operad.
\end{theorem}
\begin{proof} It suffices to provide a canonical isomorphism
\begin{equation}
\label{iso}
D\overline{\F}^\O \cong \overline{\F}^{D\O} \otimes
\overline{\mathscr H} [4-3n].
\end{equation}
To prove it, we will evaluate \eqref{iso} on an orbi-simplex $\sigma$
and establish an isomorphism, natural with respect to isomorphisms of
the corresponding graphs $\Gamma_\sigma$.

By definition,
\[
\F^{D\O}_\sigma
 = (\Gamma_\sigma^{D\O})^\vee = \bigotimes_{v\in
V(\Gamma_\sigma)}{D\O}((H(v)))^\vee
= \bigotimes_{v\in V(\Gamma_\sigma)} \bigoplus_{\substack{\text{unrooted}\\
\text{$H(v)$-trees} T_v}} \left( T_v^{\s \O[-1]} \right)^{* \vee}.
\]
Note that a graph $\Gamma_\sigma$ with each vertex $v$ decorated by a
tree $T_v$ whose leaves are labeled by the set $H(v)$ of half-edges
emanating from $v$ is literally the same as a graph $\Gamma_\tau$ with
a collection of subtrees, such that contracting each of these subtrees
returns the graph $\Gamma_\sigma$. We will call such graph
$\Gamma_\tau$ a \emph{vertex expansion} of $\Gamma_\sigma$. Moreover,
$\s \O[-1]$-decorations on the trees $T_v$ will obviously result in
$\s \O[-1]$-decorations on the graphs $\Gamma_\tau$. Thus, we see that
\[
\F^{D\O}_\sigma = \bigoplus_{\substack{\text{vertex expansions}\\
    \text{$\Gamma_\tau$ of $\Gamma_\sigma$}}} \left(\Gamma_\tau^{\s
    \O[-1]}\right)^{* \vee}.
\]
Now let us identify $\left(\Gamma^{\s \O[-1]}\right)^{*\vee}$:
\[
\left( \Gamma^{\s \O[-1]} \right)^{*\vee} = \bigotimes_{v\in
  V(\Gamma)} \Det(H(v)) \otimes \O((H(v)))^{*\vee} [3] \cong
(\Gamma^\O)^{* \vee} \otimes \bigotimes_{v\in V(\Gamma)}
\Det(H(v))[3].
\]
Leaving the factor $\left(\Gamma^\O\right)^{*\vee}$ out for the time
being, let us deal with orientations. We have
\begin{eqnarray*}
  \bigotimes_{v\in V(\Gamma)} \Det(H(v))[3]
  &
  \cong & \Det^{-3} (V(\Gamma)) \otimes \bigotimes_{v\in V(\Gamma)} \Det(H(v))\\
  &\cong& \Det^{-1} (V(\Gamma)) [2 v(\Gamma)]  \otimes \bigotimes_{v\in V(\Gamma)}
  \Det(H(v)),
\end{eqnarray*}
where $v(\Gamma) = \abs{V(\Gamma)}$. Note that the set $\coprod_{v \in
  V(\Gamma)} H(v)$ is naturally isomorphic to the set $\coprod_{e \in
  E(\Gamma)} H(e)$, where $H(e)$ is the set of (two) half-edges making
up an edge $e$, as both sets count the set of half-edges $H(\Gamma)$
of the graph, the former by grouping the set of half-edges by
vertices, the latter by edges. By passing to determinants, we obtain
\begin{equation}\label{determ}
\bigotimes_{v\in V(\Gamma)} \Det(H(v)) [3] \cong \Det^{-1}(V(\Gamma)) [2
v(\Gamma)] \otimes \bigotimes_{e\in E(\Gamma)} \Det(H(e)).
\end{equation}
Note that the exact sequence
\[
0 \to H_1(\Gamma) \to C_1(\Gamma) \to C_0(\Gamma) \to H_0(\Gamma) \to 0
\] 
yields a canonical isomorphism
\[
\Det H_0(\Gamma) \otimes \Det^{-1} H_1 (\Gamma) \cong C_0(\Gamma) \otimes
\Det^{-1} C_1 (\Gamma).
\]
Further, we have the following natural isomorphisms:
\begin{eqnarray*}
\Det C_0(\Gamma) &\cong &\Det V(\Gamma),\\
\Det C_1(\Gamma) &\cong &\bigotimes_{e \in E(\Gamma)} \Det H(e) [1]\\
& \cong &\Det^{-1} E(\Gamma) \otimes \bigotimes_{e \in E(\Gamma)} \Det H(e),\\
\Det(H_0(\Gamma)) & \cong & k[-1].
\end{eqnarray*}
We conclude that the last expression in (\ref{determ}) is isomorphic to 
\begin{align*}
  &\Det(E(\Gamma)) \otimes \Det (H_1(\Gamma))
  \otimes \Det^{-1} (H_0(\Gamma)) [2 v(\Gamma)]\\
  \cong& \Det^{-1}(E(\Gamma)) \otimes
  \Det(H_1(\Gamma)) \otimes \Det^{-1} (H_0(\Gamma)) [2 (v(\Gamma)-e(\Gamma))]\\
  \cong& \Ori^{-1}(\Gamma) [4 - 3n],
\end{align*}
which implies
\[
\F^{D\O}_\sigma \cong \bigoplus_{\text{vertex expansions $\Gamma_\tau$ of
    $\Gamma_\sigma$}} \left( \Gamma_\tau^\O \right)^{* \vee} \otimes
\Ori^{-1} (\Gamma_\tau) [4-3n].
\]
The differential on the complex $\F^{D\O}_\sigma$ is the sum of the
internal differential on $\O$ and the summation over all contractions
of a given graph $\Gamma_\tau$ along the edges arising in the vertex
expansions of $\Gamma_\sigma$ of the corresponding operad
compositions. This is similar to the differential in the graph complex $C_\bullet^{\Gamma \O}$, with
the same effect on the orientation factor, except that the resulting
grading is now cohomological.

Now let us turn to the left-hand side of \eqref{iso}. According to
Proposition \ref{GKa2} and Remark \ref{orbisimplicial}, the object $D
\F^\O_\sigma$ is represented by the complex
\begin{equation*}
  D \F^{\O}_\sigma \cong \bigoplus_{\tau\supset\sigma}
  (\F^\O_\tau \otimes \Det(\tau)[1])^* 
  = \bigoplus_{\text{vertex expansions
      $\Gamma_\tau$ of $\Gamma_\sigma$}}
  \left( \Gamma_\tau^\O \right)^{* \vee}
  \otimes \Det^{-1}(E(\Gamma_\tau)) [-1],
\end{equation*}
with the same differential as for the twisted graph complex
$\widetilde{C}_\bullet^{\Gamma \O}$. This immediately implies
\eqref{iso}.
\end{proof}
\begin{rem}
Note that the identification of $\F^{D\O}$ in the beginning of the
proof of Theorem~\ref{main2} shows that the complex computing the
cohomology of $\F^{D\O}$ is a \emph{graph complex decorated by a
decorated tree complex}, all complexes being \emph{cochain} complexes. That
complex can obviously be identified with a decorated $($cochain$)$
graph complex, and the rest of the proof of Theorem~\ref{main2}
expresses the resulting decoration through the $\O$-decoration.
\end{rem}

The following corollary describes the dualizing sheaf on $Y_n$; note
that it is concentrated in a single degree as a direct consequence of
the fact that the operad $\Comm$ is Koszul.
\begin{cor}
\label{verdier1}
The dualizing sheaf on $Y_n$ is isomorphic to ${\mathscr
  F}^{\Lie}\otimes {\mathscr H} [4 - 3n]$.
\end{cor}
\noproof

Let $\tilde{k}$ denote the one-dimensional $\Out(F_n)$-module
corresponding to the local system $\mathscr{H}$, concentrated in
degree zero; an element in $\Out(F_n)$ acts on $\tilde{k}$ as
multiplication by $1$ or $-1$, equal to the determinant of the linear
map induced on $H_1 (\Gamma)$.  Then we have the following result
which was formulated by Kontsevich in \cite{Kon} and given a different
proof in \cite{CV}.
\begin{cor}
\label{shift1}
There are the following isomorphisms of graded $k$-vector spaces:
\begin{enumerate}
\item
$H_\bullet(\Out(F_n),k)\cong H_{\Gamma \Lie}^{3n-4 - \bullet} (n)$;
\item
$H_\bullet(\Out(F_n),\tilde{k})\cong \widetilde{H}_{\Gamma
    \Lie}^{3n-4 - \bullet} (n)$.
\end{enumerate}
\end{cor}
\begin{proof}
  As usual, we limit ourselves with proving the first statement. We
  have
\begin{align*}
  H_\bullet(\Out(F_n),k)&  \cong H_\bullet(Y_n, k)\\
  &  \cong [H^\bullet(Y_n, k)]^\vee\\
  &  \cong [H^\bullet(Y_n, \F^{\Comm})]^\vee\\
  &\cong H_c^\bullet(Y_n, D\mathscr F^{\Comm})^{* \vee}\\
  &\cong H_c^\bullet(Y_n, \mathscr F^{D\Comm}\otimes \mathscr H [4-3n])^{* \vee}\\
  &\cong H_c^{\bullet+4 - 3n}(Y_n, \mathscr F^{D\Comm}\otimes \mathscr H)^{* \vee}\\
  &\cong H_c^{3n-4 - \bullet}(Y_n,\mathscr F^{\Lie}\otimes \mathscr H)\\
  &\cong H^{3n-4 - \bullet}_{\Gamma \Lie},
\end{align*}
as required. Note that we used the fact that $D = D^{-1}$ where $D$ is
the functor of taking the Verdier dual.
\end{proof}

\begin{rem}
  Compare this to a more straightforward computation to get a relation
  between the cohomology of $Y_n$ with compact supports and
  commutative graph cohomology:
\[
H_c^\bullet(Y_n, k) = H_c^\bullet(Y_n, \F^{\Comm}) =
\widetilde{H}^{\bullet}_{\Gamma \Comm} (n).
\]
\end{rem}

\section{Ribbon graphs}

The theory developed in the previous section has an analogue for
non-$\Sigma$ operads and ribbon graph complexes. Recall that a ribbon
graph is an $\Ass$-decorated graph; this is equivalent to having a
cyclic ordering on the set of half-edges around each vertex. Given a
ribbon graph $\Gamma$, there is a canonical way of producing a
compact, oriented surface with boundary $S(\Gamma)$ of which the graph
$\Gamma$ is a deformation retract. In this way one attaches to a
ribbon graph two invariants: the genus $g \ge 0$ and the number $n
\ge 1$ of boundary components of the corresponding surface, $2-2g-n <
0$. An isomorphism between two ribbon graphs is an isomorphism
preserving the cyclic ordering around each vertex. We will not specify
whether the boundary components should be fixed (not necessarily
point-wise) under an isomorphism or allowed to be permuted; both
versions admit completely parallel treatments.

The \emph{mapping class group} $\Gamma_{g,n}$ is the group of isotopy
classes of orientation-preserving diffeomorphisms of an oriented
surface of genus $g$ with $n$ boundary components. Again, we shall be
ambiguous whether or not $\Gamma_{g,n}$ permutes the boundary components of a
surface.

Now let $\O$ be a cyclic (chain) $k$-operad with $\O(1) = k$ without
the action of the symmetric group, a so-called non-$\Sigma$ operad.
We introduce the notions of a ribbon $\O$-graph complex
$C_\bullet^{\Rib \O}$, its (cochain) dual $C^\bullet_{\Rib \O}$, as
well as the twisted versions $\widetilde{C}_\bullet^{\Rib \O}$ and
$\widetilde{C}^\bullet_{\Rib \O}$ in precisely the same way as in the
previous section. For two quasi-isomorphic operads, the corresponding
ribbon graph complexes will be quasi-isomorphic. The subcomplex in
$C_\bullet^{\Rib \O}$ consisting of ribbon graphs with fixed $g$ and
$n$ will be denoted by $C_\bullet^{\Rib \O}(g,n)$ ($C^\bullet_{\Rib
  \O}(g,n)$ for the cohomological version). It is easy to see that
\[
C_\bullet^{\Rib \O} \cong \bigoplus_{\substack{g\geq 0, n\geq 1 \\
    2-2g-n < 0}}C_\bullet^{\Rib \O}(g,n).
\]
The most important example of a ribbon $\O$-graph complex corresponds
to the associative non-$\Sigma$ operad $\T(m)=k$ in all degrees $m \ge
1$.  In this case we have an isomorphism $C_\bullet^{\Rib \T}(g,n)
\cong C_\bullet^{\Gamma \Ass}(g,n)$ and similarly for the twisted
versions, cf.\ the discussion at the end of Section~\ref{graphs}.

We will now introduce the space of metric ribbon graphs $\M_{g,n}$. A
point $\Gamma$ in $\M{g,n}$ is an isomorphism class of ribbon graphs
of genus $g$ with $n$ boundary components such that each edge $e\in
E(\Gamma)$ is supplied with length $l(e)>0$; we require that
$\sum_{e\in E(\Gamma)}l(e)=1$. The space $\M_{g,n}$ naturally
compactifies to a simplicial orbi-complex $\Mbar_{g,n}$. Namely, one
introduces a contractible space of \emph{marked ribbon graphs}
$\A_{g,n}$, also known as the \emph{arc complex}, on which
$\Gamma_{g,n}$ acts properly discontinuously. The space $\A_{g,n}$
admits a simplicial closure $\overline{\A}_{g,n}$ (also known as the \emph{arc complex}) and we denote by
$\Mbar_{g,n}$ the quotient of $\overline{\A}_{g,n}$ by $\Gamma_{g,n}$.

\begin{rem}
  The space $\M_{g,n}$ is known to be homeomorphic to the Cartesian
  product of the open standard simplex $\Delta^{n-1}$ and the moduli
  space of Riemann surfaces of genus $g$ with $n$ labeled punctures or
  its quotient by the diagonal action of the symmetric group $S_n$,
  depending on whether we allow graph isomorphisms permuting the
  boundary components. A version of the simplicial compactification of
  $\M_{g,n}$ was constructed by Kontsevich in \cite{Kon'} in
  connection with his proof of the Witten conjecture.  The reader is
  referred to the papers by Looijenga \cite{Loo}, Zvonkine \cite{Zv},
  and Mondello \cite{Mon} for details and comparison of $\Mbar_{g,n}$
  to the Deligne-Mumford compactification.
\end{rem}

Definitions of the sheaves $\F^\O$, $\overline{\F}^\O$, and $\mathscr
H$ on $\M_{g,n}$ and $\Mbar_{g,n}$ transfer verbatim from the
corresponding definitions in the previous section. We can now
formulate analogues of the main results from Section~\ref{metric}.
They are proved in precisely the same way as Theorems \ref{main1} and
\ref{main2}.
\begin{theorem}
\label{main3}
There are canonical isomorphisms of $k$-vector spaces:
\begin{enumerate}
\item $\widetilde{H}^\bullet_{\Rib\O}(g,n)\cong
  H^{\bullet}_c(\M_{g,n}, \F^{\O})$;
\item ${H}^\bullet_{\Rib\O}(g,n) \cong H^{\bullet}_c(\M_{g,n},
  \F^{\O}\otimes \mathscr H)$.
\end{enumerate}
\end{theorem}
\noproof

Note the following explicit relation between the cohomology of $Y_r$
and $\M_{g,n}$.
\begin{cor}
\[
H_c^\bullet (Y_r, \F^{\Ass}) \cong \bigoplus_{\substack{g \ge 0, n \ge 1\\
    2 - 2g - n = 1 - r}} H_c^\bullet (\M_{g,n}, k) \qquad \text{for $r >
  1$}.
\]
\end{cor}
\noproof

\begin{theorem}
\label{main4}
There is a canonical isomorphism in the derived category of sheaves on
$\M_{g,n}$:
\[
D{\F^\O}\cong {\F}^{D\O}\otimes {\mathscr H} [7-6g-3n].
\]
Here $D\mathscr O$ is the dg-dual non-$\Sigma$ operad to $\mathscr O$,
defined for non-$\Sigma$ operads the same way as in $\eqref{dg-dual}$,
except that the trees must be planar.
\end{theorem}
\noproof
The analogue of Corollary \ref{shift1} reads as follows.
\begin{cor}
\label{shift2}
There are the following isomorphisms of graded $k$-vector spaces:
\begin{enumerate}
\item $H_\bullet(\M_{g,n},k) \cong H_{\Rib \T}^{6g+3n-7-\bullet} (g,n)
  \cong H_{\Gamma \Ass}^{6g+3n-7-\bullet} (g,n)$;
\item $H_\bullet(\M_{g,n},\tilde{k}) \cong \widetilde{H}_{\Rib
    \T}^{6g+3n-7-\bullet} (g,n) \cong \widetilde{H}_{\Gamma \Ass}^{6g+
    3n-7 - \bullet} (g,n)$.
\end{enumerate}
\end{cor}
\noproof The following corollary describes the dualizing sheaf on
$\M_{g,n}$.
\begin{cor}
  \label{verdier2}
  The dualizing sheaf on $\M_{g,n}$ is isomorphic to ${\mathscr H}
  [7-6g-3n]$.
\end{cor}
\noproof
\begin{rem}
  Note that the dualizing sheaf is locally constant in agreement with
  the well-known fact that $\M_{g,n}$ is an orbifold.  Therefore,
  Verdier duality on $\M_{g,n}$ turns into Poincar\'e-Lefschetz
  duality:
\[
H_{\bullet} (\M_{g,n}, k) = H^{6g+3n-7 - \bullet}_c (\M_{g,n},
\tilde{k}).
\]
When $\M_{g,n}$ stands for the moduli space of ribbon graphs with
\emph{labeled} boundary components, it will be an orientable orbifold,
in which case $\tilde k \cong k$, but taking its quotient by the
symmetric group permuting the boundary components to get the other
version of $\M_{g,n}$ will destroy orientability, and the dualizing
sheaf will no longer be constant.
\end{rem}

\end{document}